\documentclass[proceedings,submission%
]{dmtcs}

\usepackage[latin1]{inputenc}
\usepackage{subfigure,wrapfig}

\usepackage[round]{natbib} 
 
\usepackage{amsmath,amssymb,amsfonts,latexsym, array, float}
\usepackage{tikz}
\usetikzlibrary{shadows}

\usepackage{graphicx}
\newcommand{\params}[1]{$($#1\rotatebox[origin=c]{180}{$($}}

\newtheorem{defi}{Definition}[section]
\newtheorem{theo}[defi]{Theorem}

\newtheorem{example}[defi]{Example}
\usepackage[ruled,vlined]{algorithm2e}

\usepackage{xspace}

\newcommand{\C}{\ensuremath{\mathcal{C}}\xspace}

\newcommand{\Sys}{\ensuremath{\mathcal{E}}\xspace}
\def\cal{\mathcal}

\newcommand{\D}{\ensuremath{{\mathcal D}}}

\newcommand{\s}{\ensuremath{{\mathcal S}}}

\newcommand{\Bstar}{\ensuremath{{B^\star}}}
\newcommand{\wc}{\ensuremath{{\hat{\mathcal{C}}}}}
\newcommand{\Wc}{\ensuremath{{\hat{\mathcal{C}}}}}
\newcommand{\rs}{restriction\xspace}
\newcommand{\rss}{restrictions\xspace}
\newcommand{\rt}{restriction term\xspace}
\newcommand{\rts}{restriction terms\xspace}
\begin{document}

\title{Combinatorial specification of permutation classes\thanks{This work was completed with the support of the ANR project MAGNUM number 2010\_BLAN\_0204.}}

\author[ F. Bassino \and 
        M. Bouvel \and 
        A. Pierrot \and 
        C. Pivoteau \and 
        D. Rossin]{
        Fr\'ed\'erique Bassino\addressmark{1} \and 
        Mathilde Bouvel\addressmark{2} \and 
        Adeline Pierrot\addressmark{3} \and 
        Carine Pivoteau\addressmark{4} \and 
        Dominique Rossin\addressmark{5} 
}


\address{
\addressmark{1}Universit\'e Paris 13, LIPN (CNRS UMR 7030), Villetaneuse, France.\\
\addressmark{2}Universit\'e de Bordeaux, LaBRI (CNRS UMR 5800), Talence, France.\\
\addressmark{3}Universit\'e Paris Diderot, LIAFA (CNRS UMR 7089), Paris, France.\\
\addressmark{4}Universit\'e Paris-Est, LIGM (CNRS UMR 8049), Marne-la-Vall\'ee, France.\\
\addressmark{5}\'Ecole Polytechnique, LIX (CNRS UMR 7161), Palaiseau, France.
}


\maketitle
\begin{abstract}
\paragraph{Abstract.}
This article presents a methodology that automatically derives a 
combinatorial specification for the permutation class $\mathcal{C} 
= Av(B)$, given its basis $B$ of excluded patterns and the 
set of simple permutations in $\mathcal{C}$, when these sets are both 
finite. 
This is achieved considering both pattern avoidance and pattern 
containment constraints in permutations.
The obtained specification yields a system of equations satisfied 
by the generating function of $\mathcal{C}$, this system being always positive
and algebraic. It also yields a uniform random sampler of permutations 
in $\mathcal{C}$. The method presented
is fully algorithmic.

\paragraph{R\'esum\'e.}
Cet article pr\'esente une m\'ethodologie qui calcule automatiquement une
sp\'ecification combinatoire pour la classe de permutations
$\mathcal{C} = Av(B)$, \'etant donn\'es une base $B$ de motifs interdits
et l'ensemble des permutations simples de $\mathcal{C}$, lorsque ces deux
ensembles sont finis. 
Ce r\'esultat est obtenu en consid\'erant \`a la fois des contraintes de 
motifs interdits et de motifs obligatoires dans les permutations.
La sp\'ecification obtenue donne un syst\`eme d'\'equations satisfait
par la s\'erie g\'en\'eratrice de la classe $\mathcal{C}$, syst\`eme qui est 
toujours positif et alg\'ebrique.
Elle fournit aussi un g\'en\'erateur al\'eatoire uniforme de permutations 
dans $\mathcal{C}$. La m\'ethode pr\'esent\'ee est compl\`etement
algorithmique.
\end{abstract}

\section{Introduction}
\label{sec:intro}

Initiated by \cite{Knuth73} almost forty years ago, the study of permutation 
classes has since received a lot of attention, mostly with respect to enumerative questions (see
\cite{Bou03,Eli04,KiMa03} and their references among many others).
Most articles are focused on a given class $\mathcal{C} = Av(B)$ where
the basis $B$ of excluded patterns characterizing $\mathcal{C}$ is
finite, explicit, and in most cases contains only patterns of size $3$
or $4$.  Recently, some results of a rather different nature have been
obtained, and have in common that they describe general properties of
permutation classes -- see
\cite{AA05,ALR05,BBPR09,BBPR11,BHV06a,BRV06,Vat05} for example.
Our work falls into this new line of research.

Our goal in this article is to provide a general algorithmic method to
obtain a combinatorial specification for any permutation class
$\mathcal{C}$ from its basis $B$ and the set $\mathcal{S}_{\mathcal{C}}$
of simple permutations in $\mathcal{C}$, and assuming these two sets are
finite.
Notice that by previous works to be detailed in Section 
\ref{sec:previous_algos}, it is enough to know the finite basis $B$ of the
class to decide whether the set $\mathcal{S}_{\mathcal{C}}$ is finite 
and (in the affirmative) to compute $\mathcal{S}_{\mathcal{C}}$.

By \emph{combinatorial specification} of a class (see \cite{FlSe09}), we
mean an unambiguous system of combinatorial equations that describe
recursively the permutations of $\mathcal{C}$ using only combinatorial constructors
(disjoint union, cartesian product, sequence, \ldots) and permutations
of size $1$. 
Notice the major difference with the results of \cite{AA05}: our 
specifications are unambiguous, whereas \cite{AA05} obtain combinatorial 
systems of equations characterizing permutations classes that are ambiguous 
in general.

We believe that our purpose of obtaining algorithmically combinatorial 
specifications of permutation classes is of interest \emph{per se} but
also because it then allows to obtain by routine algorithms a system
of equations satisfied by the generating function of $\mathcal{C}$ and
a Boltzmann uniform random sampler of permutations in $\mathcal{C}$,
using the methods of \cite{FlSe09} and \cite{DuFlLoSc04} respectively.

The paper is organized as follows.  
Section~\ref{sec:PermClasses} proceeds with some background on 
permutation classes, simple permutations and substitution decomposition, and
Section~\ref{sec:previous_algos} sets the algorithmic context of our study.
Section~\ref{sec:ambiguous} then explains how to obtain a system of
combinatorial equations describing $\mathcal{C}$ from the set
of simple permutations in $\mathcal{C}$, that we assume to be finite.
The system so obtained may
be ambiguous and Section~\ref{sec:disambiguation} describes a
disambiguation algorithm to obtain a combinatorial specification for
$\mathcal{C}$. 
The most important idea of this disambiguation procedure is to 
transform ambiguous unions into disjoint unions of terms that involve both 
pattern avoidance and pattern containment constraints. This somehow
allows to interpret on the combinatorial objects themselves the result of 
applying the inclusion-exclusion on their generating functions.
Finally, Section~\ref{sec:ccl} concludes the whole
algorithmic process by explaining how this specification can be
plugged into the general methodologies of \cite{FlSe09} and
\cite{DuFlLoSc04} to obtain a system of equations satisfied by the
generating function of $\mathcal{C}$ and a Boltzmann uniform random
sampler of permutations in $\mathcal{C}$. 
We also give a number of perspectives opened by our algorithm.

\section{Permutation classes and simple permutations}
\label{sec:PermClasses}

\subsection{Permutation patterns and permutation classes}

A permutation $\sigma= \sigma_1 \sigma_2 \ldots \sigma_n$ of size
$|\sigma| = n$ is a bijective map from $\{1,\ldots ,n\}$ to itself,
each $\sigma_i$ denoting the image of $i$ under $\sigma$.  A
permutation $\pi = \pi_1 \pi_2 \ldots \pi_k$ is a \emph{pattern} of a
permutation ${\sigma = \sigma_1 \sigma_2 \ldots \sigma_n}$ (denoted $\pi
\preceq \sigma$) if and only if $k\leq n$ and there exist integers $1
\leq i_1 < i_2 < \ldots < i_k \leq n$ such that $\sigma_{i_1}\ldots
\sigma_{i_k}$ is order-isomorphic to $\pi$, \emph{i.e.} such that
$\sigma_{i_{\ell}} < \sigma_{i_m}$ whenever $\pi_{\ell} < \pi_{m}$.  A
permutation~$\sigma$ that does not contain $\pi$ as a pattern is said
to {\em avoid} $\pi$.  For example the permutation $\sigma=316452$
contains~$\pi = 2431$ as a pattern, whose occurrences are $3642$ and
$3652$. But $\sigma$ avoids the pattern $2413$ as none of its
subsequences of length $4$ is order-isomorphic to $2413$.

The pattern containment relation $\preceq$ is a partial order on
permutations, and a {\em permutation class} $\mathcal{C}$ is a downset
under this order: for any $\sigma \in \mathcal{C}$, if $\pi \preceq
\sigma$, then we also have $\pi \in \mathcal{C}$.  For every set $B$,
the set $Av(B)$ of permutations avoiding any pattern of $B$ is a
class. Furthermore every class $\mathcal{C}$ can be rewritten as
$\mathcal{C} = Av(B)$ for a unique antichain $B$ ({\em i.e.,} a unique
set of pairwise incomparable elements) called the {\it basis} of
$\mathcal{C}$. 
The basis of a class $\mathcal{C}$ may be finite or infinite; it is
described as the set of permutations that do not belong to $\mathcal{C}$ and
that are minimal in the sense of $\preceq$ for this criterion.

In the following, we only consider classes whose basis $B$ is given
explicitly, and is finite. This does not cover the whole range of
permutation classes, but it is a reasonable assumption when dealing
with \emph{algorithms} on permutation classes, that take a finite
description of a permutation class as input.
Moreover, as proved by \cite{AA05}, it is necessary that $B$ is finite as soon 
as the set $\mathcal{S}_{\mathcal{C}}$ of simple permutations in $\mathcal{C}
 = Av(B)$ is finite. Consequently the assumption of the 
finiteness of $B$ is not a restriction when working on permutation classes 
such that $\mathcal{S}_{\mathcal{C}}$ is finite, which is the context of 
our study.

\subsection{Simple permutations and substitution decomposition of permutations}
\label{sec:substitution}

An \emph{interval} (or {\em block}) of a permutation $\sigma$ of size
$n$ is a subset $\{i,\ldots ,(i+\ell-1)\}$ of consecutive integers of
$\{1,\ldots ,n\}$ whose images by $\sigma$ also form an interval of
$\{1,\ldots ,n\}$. The integer $\ell$ is called the \emph{size} of the
interval.  A permutation $\sigma$ is \emph{simple} when it is of size
at least $4$ and it contains no interval, except the trivial ones:
those of size $1$ (the singletons) or of size $n$ ($\sigma$
itself). The permutations $1$, $12$ and $21$ also have only trivial
intervals, nevertheless they are \emph{not} considered to be simple
here. Moreover no permutation of size $3$ has only trivial intervals.
For a detailed study of simple permutations, in particular from an
enumerative point of view, we refer the reader to
\cite{AA05,AAK03,Bri08}. 

Let $\sigma$ be a permutation of size $n$ and $\pi^{1},\ldots,
\pi^{n}$ be $n$ permutations of size $p_1, \ldots, p_n$
respectively. Define the \emph{substitution} $\sigma[\pi^{1}, \pi^{2}
,\ldots, \pi^{n}]$ of $\pi^{1},\pi^{2} , \ldots, \pi^{n}$ in
$\sigma$ 
to be the permutation of size $p_1 + \ldots + p_n$ obtained by
concatenation of $n$ sequences of integers $S^1, \ldots , S^n$ from
left to right, such that for every $i,j$, the integers of $S^i$ form
an interval, are ordered in a sequence order-isomorphic to $\pi^{i}$,
and $S^i$ consists of integers smaller than $S^j$ if and only if
$\sigma_i < \sigma_j$.  For instance, the substitution $ 1\, 3\, 2
[2\, 1, 1\, 3\, 2, 1]$ gives the permutation $ 2\, 1\, \, 4\, 6\, 5\,
\, 3$.  We say that a permutation $\pi$ is \emph{$12$-indecomposable}
(resp. \emph{$21$-indecomposable}) if it cannot be written as
$12[\pi^{1},\pi^{2}]$ (resp. $21[\pi^{1},\pi^{2}]$), for any
permutations $\pi^{1}$ and $\pi^{2}$.

Simple permutations allow to describe all permutations through their \emph{substitution decomposition}.
\begin{theo}[\cite{AA05}]
Every permutation $\pi$ of size $n$ with
  $n \geq 2$ can be uniquely decomposed as follows, $12$ (resp. $21$, $\sigma$) being called the \emph{root} of $\pi$:
  
\vspace{-5.5pt}
\begin{itemize}\setlength{\itemsep}{-0.8pt}
\item $12[\pi^{1},\pi^{2}]$, with $\pi^{1}$ $12$-indecomposable,
\item $21[\pi^{1},\pi^{2}]$, with $\pi^{1}$ $21$-indecomposable,
\item $\sigma[\pi^{1},\pi^{2},\ldots,\pi^{k}]$, with $\sigma$ a simple permutation of size $k$.
\end{itemize}
\label{thm:decomp_perm_AA05}
\end{theo}
\vspace{-9pt}

To account for the first two items of
Theorem~\ref{thm:decomp_perm_AA05} in later discussions, we
furthermore introduce the following notations: For any set $\C$ of
permutations, $\C^+$ (resp. $\C^-$) denotes the set of permutations of
$\C$ that are $12$-indecomposable (resp. $21$-indecomposable). Notice
that even when $\C$ is a permutation class, this is not the case for
$\C^+$ and $\C^-$ in general.

Theorem~\ref{thm:decomp_perm_AA05} provides the first step in the
decomposition of a permutation $\pi$. To obtain its full
decomposition, we can recursively decompose the permutations $\pi^{i}$
in the same fashion, until we reach permutations of size $1$. This
recursive decomposition can naturally be represented by a tree, that
is called the substitution decomposition tree (or {\em decomposition
  tree} for short) of $\pi$.  Each internal node of the tree is
labeled by $12,21$ or by a simple permutation and the leaves represent
permutation $1$.
Notice that in decomposition trees, the left child of a 
node labeled $12$ (resp. $21$) is never labeled by $12$ (resp. $21$),
since $\pi^{1}$ is $12$-indecomposable (resp. $21$-indecomposable) in
the first (resp. second) item of Theorem~\ref{thm:decomp_perm_AA05}.

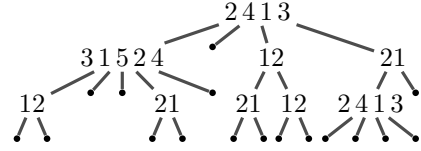
\begin{wrapfigure}[8]{r}{60mm}
        \vspace{-3.5mm}
\begin{tikzpicture}[
 scale=.6,
    level/.style={sibling distance=20mm/#1},
   edge from parent/.style={very thick,draw=black!70},
    simple/.style={rectangle, draw=none, rounded corners=1mm, fill=white, text centered, text=black,anchor=north,inner sep=2pt},
    linear/.style={rectangle, draw=none, fill=white, text centered, anchor=north, text=black,inner sep=2pt},
    every node/.style={circle, draw=none, fill=black, text centered, anchor=north, text=white,inner sep=0},
    level distance=7mm
]
\node[simple] {$2\,4\,1\,3$}
	child { node[simple] {$3\,1\,5\,2\,4$}
		child {node[linear] {$12$}
			child { node { ~ }}
			child { node { ~ }}
		}
		child[sibling distance=7mm] { node { ~ }}
		child { node { ~ }}
		child{node[linear] {$21$}
			child { node { ~ }}	
			child { node { ~ }}	
		}
		child { node { ~ }}
	}
	child {node { ~ }}
	child[sibling distance=6mm] {node[linear] {$12$}
		child {node[linear] {$21$}
        		child { node { ~ }}
        		child { node { ~ }}
        	}
		child {node[linear] {$12$}
        		child { node { ~ }}
        		child { node { ~ }}
        	}
	}
	child {node[linear] {$21$}
	   child {node[simple] {$2\,4\,1\,3$}
			child { node { ~ }}
			child { node { ~ }}
			child { node { ~ }}
			child { node { ~ }}
		}
		child {node { ~ }}
};
\end{tikzpicture}
\caption{{\small Decomposition tree of $\pi$ (from Ex.~\ref{ex:decomp}).}}\label{fig:tree}
\end{wrapfigure}

\begin{example}\label{ex:decomp}
  The permutation $\pi = 8\ 9\ 5\ 11\ 7\ 6\ 10\ 17\ 2\ 1\ 3\ 4\ 14\
  16\ 13\ 15\ 12$ is recursively decomposed as $\pi =
  2413[4517326,1,2134,35241] =
  2413[31524[12[1,1],1,1,21[1,1],1]],1,12[21[1,1],12[1,1]],\\
  21[2413[1,1,1,1],1]] $ and its decomposition tree is given in
  Figure~\ref{fig:tree}.
\end{example}

The \emph{substitution closure} $\Wc$ of a permutation
class\footnote{that contains permutations $12$ and $21$. We will
  assume so in the rest of this article to avoid trivial cases.}
$\mathcal{C}$ is defined as the set of permutations whose
decomposition trees have internal nodes labeled by either $12, 21$ or
a simple permutation of $\mathcal{C}$. Notice that $\mathcal{C}$ and
$\Wc$ therefore contain the same simple permutations. Obviously, for
any class $\mathcal{C}$, we have $\mathcal{C} \subseteq \Wc$. When the
equality holds, the class $\mathcal{C}$ is said to be {\em
  substitution-closed} (or sometimes {\em wreath-closed}). But this is
not always the case, and the simplest example is given by $\mathcal{C}
= Av(213)$. This class contains no simple permutation hence its
substitution closure is the class of separable permutations of
\cite{BBL98}, \emph{i.e.} of permutations whose decomposition trees
have internal nodes labeled by $12$ and $21$. It is immediate to
notice that $213 \in \Wc$ whereas of course $213 \notin \mathcal{C}$.

A characterization of substitution-closed classes useful for our
purpose is given in \cite{AA05}: A class is substitution-closed if and
only if its basis contains only simple permutations.

\section{Algorithmic context of our work}
\label{sec:previous_algos}

\begin{figure}[ht]
                \includegraphics[width=\textwidth]{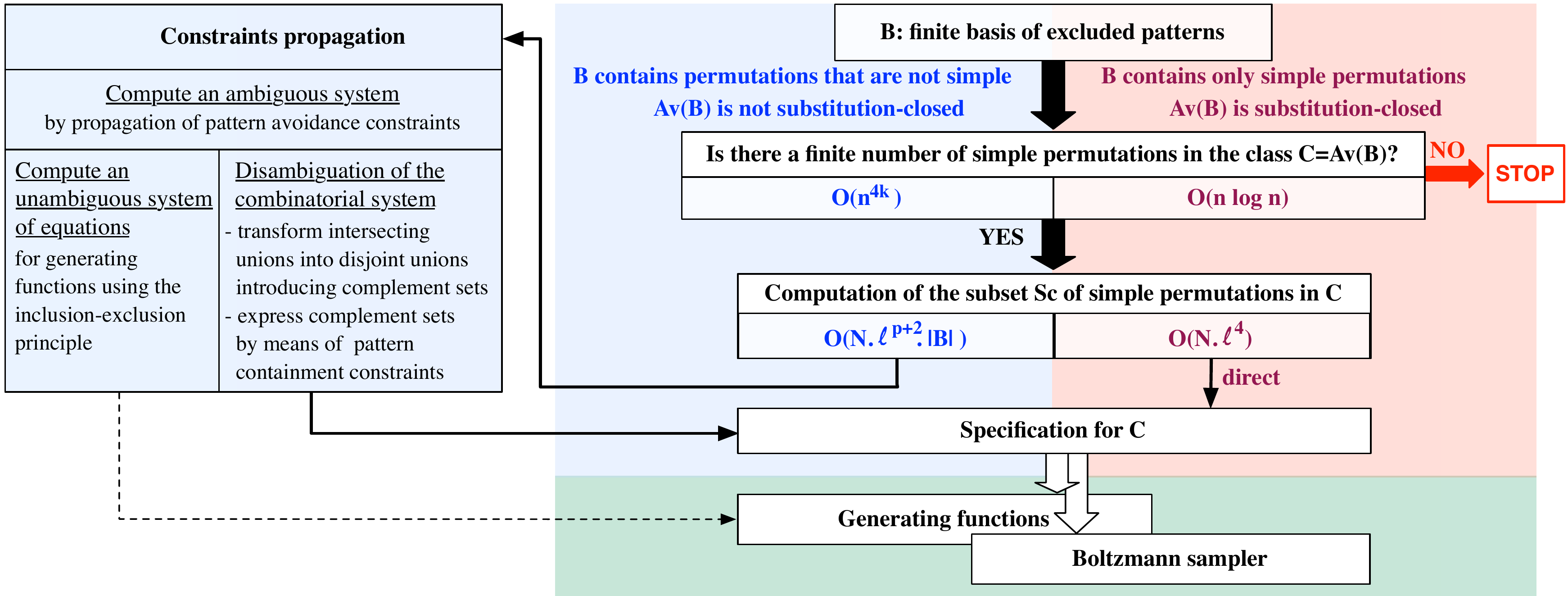}
		         \vspace{-3mm}
                \caption{Automatic process from the basis of a
                  permutation class to generating function and
                  Boltzmann sampler.}
        \label{fig:schema2}
\end{figure}

Putting together the work reported in this article and recent algorithms
from the litterature provides a full algorithmic chain starting with the
finite basis $B$ of a permutation class $\mathcal{C}$, and computing a 
specification for $\mathcal{C}$. 
The hope for such a very general algorithm is of course very tenuous,
and the algorithm we describe below will compute its output only when 
some hypothesis are satisfied, which are also tested algorithmically.
Figure~\ref{fig:schema2} summarizes the main steps of the algorithm.

The algorithms 
performing the first two steps of the algorithmic process of Figure~\ref{fig:schema2}
are as follows. 

\textbf{First step : Finite number of simple permutations} \qquad
First, we check whether $\mathcal{C} = Av(B)$ contains only a finite number of
simple permutations.  
This is achieved using algorithms of
\cite{BBPR09} when the class is substitution-closed and of
\cite{BBPR11} otherwise. The complexity of these algorithms are
respectively $\mathcal{O}(n \log n)$ and $\mathcal{O}(n^{4k})$, where
$n = \sum_{\beta \in B} |\beta|$ and $k = |B|$.

 \textbf{Second step : Computing simple permutations} \qquad
 The second step of the algorithm is the computation of the set of
 simple permutations $\mathcal{S}_{\mathcal{C}}$ contained in
 $\mathcal{C} = Av(B)$, when we know it is finite. Again, when
 $\mathcal{C}$ is substitution-closed, $\mathcal{S}_{\mathcal{C}}$ can
 be computed by an algorithm that is more efficient than in the
 general case. The two algorithms are described in \cite{PR11}, and
 their complexity depends on the output: $\mathcal{O}(N \cdot
 \ell^{p+2}\cdot |B|)$ in general and $\mathcal{O}(N \cdot \ell^{4})$ for
 substitution-closed classes, with $N = |\mathcal{S}_{\mathcal{C}}|$,
 $p = \max \{|\beta| : \beta \in B\}$ and $\ell = \max \{|\pi| : \pi
 \in \mathcal{S}_{\mathcal{C}}\}$.
 
Sections~\ref{sec:ambiguous}
and~\ref{sec:disambiguation} will then explain how to derive a
specification for $\mathcal{C}$ from $\mathcal{S}_{\mathcal{C}}$.

\section{Ambiguous combinatorial system describing $\mathcal{C}$\label{sec:ambigu}}
\label{sec:ambiguous}

We describe here an algorithm that 
takes as input the set $\mathcal{S}_{\mathcal{C}}$ of
simple permutations in a class $\mathcal{C}$ and the basis $B$ of $\mathcal{C}$, 
and that produces in output a
(possibly ambiguous) system of combinatorial equations describing the
permutations of $\mathcal{C}$ through their decomposition trees.
The main ideas are those of Theorem~10 of \cite{AA05}, but unlike this work,
we make the whole process fully algorithmic.

\subsection{The simple case of substitution-closed classes}

Recall that $\mathcal{C}$ is a substitution-closed permutation class
when {$\mathcal{C}=\Wc$}, or equivalently when the permutations in
$\mathcal{C}$ are exactly the ones whose decomposition trees have
internal nodes labeled by $12, 21$ or any simple permutation of
{$\C$}.  Then Theorem~\ref{thm:decomp_perm_AA05} directly yields the
following system $\Sys_{\wc}$:
        \begin{eqnarray}
        \wc  &=& 1\ \uplus \ 12[\wc^+, \wc]\  \uplus \ 21[\wc^-, \wc] \ \textstyle\uplus \biguplus_{\pi \in \s_\wc} \pi[\wc, \dots, \wc]\label{eqn:Wc1}  \\  
        \wc^+ &=& 1\  \uplus \ 21[\wc^-, \wc]\  \uplus\  \textstyle\biguplus_{\pi \in \s_\wc} \pi[\wc, \dots, \wc] \label{eqn:Wc2}\\        
        \wc^- &=& 1\  \uplus\  12[\wc^+, \wc]\  \uplus\  \textstyle\biguplus_{\pi \in \s_\wc} \pi[\wc, \dots, \wc]. \label{eqn:Wc3}
        \end{eqnarray}
        By uniqueness of substitution decomposition, unions are disjoint
        and so Equations~\eqref{eqn:Wc1} to \eqref{eqn:Wc3} describe
        unambiguously the substitution closure $\Wc$ of a permutation
        class $\mathcal{C}$. For a substitution-closed class (and the
        substitution closure of any class), this description gives a
        combinatorial specification. Hence, it provides an efficient
        way to compute the generating function of the class, and to
        generate uniformly at random a permutation of a given size in
        the class.

\subsection{Adding constraints for classes that are not substitution-closed\label{sec:addConstraint}}

When ${\mathcal C}$ is not substitution-closed, we
compute a new system by adding constraints to the system
obtained for \wc, as in \cite{AA05}. Denoting by $X\langle Y\rangle$ the
set of permutations of $X$ that avoid the patterns in $Y$, we
have $\C = \wc\langle\Bstar\rangle$ where $\Bstar$ is the
subset of non-simple permutations of~$B$. Noticing that $\s_\wc = \s_\C$ (by definition of $\wc$), and since 
$\C^{\varepsilon} = \wc^{\varepsilon} \langle \Bstar
\rangle$ 
for $\varepsilon \in \{~~ , +, -\}$
, Equations~\eqref{eqn:Wc1} to~\eqref{eqn:Wc3} give
        \begin{eqnarray}
          \wc \langle \Bstar \rangle &=& 1\ \uplus \ 12[\wc^+, \wc]\langle \Bstar \rangle\  \uplus \ 21[\wc^-, \wc]\langle \Bstar \rangle \ \uplus \textstyle\biguplus_{\pi \in \s_\C} \pi[\wc, \dots, \wc]\langle \Bstar \rangle \label{eqn:const1}\\       
          \wc^+ \langle \Bstar \rangle &=& 1\  \uplus \ 21[\wc^-, \wc]\langle \Bstar \rangle\  \uplus\  \textstyle\biguplus_{\pi \in \s_\C} \pi[\wc, \dots, \wc]\langle \Bstar \rangle \label{eqn:const2}\\        
          \wc^- \langle \Bstar \rangle &=& 1\  \uplus\  12[\wc^+, \wc]\langle \Bstar \rangle\  \uplus\  \textstyle\biguplus_{\pi \in \s_\C} \pi[\wc, \dots, \wc]\langle \Bstar \rangle, \label{eqn:const3}
        \end{eqnarray}
all these unions being disjoint. 
This specification is not complete, since sets of the form $\pi[\wc, \dots, \wc]\langle \Bstar \rangle$ are not immediately described  from $\wc \langle \Bstar\rangle$. 
Theorem 10 of \cite{AA05} explains how sets such as $ \pi[\wc, \dots, \wc]\langle \Bstar \rangle$ can be expressed as union of smaller sets:
$$\pi[\wc, \dots, \wc]\langle \Bstar \rangle = \textstyle\bigcup_{i=1}^{k} \pi[\wc\langle E_{i,1} \rangle,\wc\langle E_{i,2} \rangle,\ldots,\wc\langle E_{i,k} \rangle]$$
where $E_{i,j}$ are sets of permutations which are patterns of some
permutations of $\Bstar$. 
This introduces sets of the form $\wc\langle E_{i,j} \rangle$ on the right-hand side of an equation of the system that do not appear on the left-hand side of any equation. We will call such sets \emph{right-only} sets.
Taking $E_{i,j}$ instead of $\Bstar$ in Equations~\eqref{eqn:const1} to~\eqref{eqn:const3}, we can recursively compute these right-only sets by introducing new equations in the system. 
This process terminates since there exists only a finite number of sets of patterns of elements of $\Bstar$ (as $B$ is finite). Let us introduce some definitions to describe these sets $E_{i,j}$.

A \emph{generalized substitution} $\sigma\{\pi^{1}, \pi^{2} ,\ldots,
\pi^{n}\}$ is defined as a substitution (see
p.\pageref{thm:decomp_perm_AA05}) with the particularity that any
$\pi^i$ may be the empty permutation (denoted by~$0$). Specifically
$\sigma[\pi^{1}, \pi^{2} ,\ldots, \pi^{n}]$ necessarily contains
$\sigma$ whereas $\sigma\{\pi^{1}, \pi^{2} ,\ldots, \pi^{n}\}$ may
avoid $\sigma$.  For instance, $ 1\, 3\, 2 \{2\, 1, 0, 1\}= 2\, 1\,
3\in Av(132)$.

An {\em embedding of $\gamma$ in~${\pi=\pi_1\dots\pi_n}$} is a
map $\alpha$ from $\{1, \ldots, n\}$ to the set of (possibly empty) 
blocks\footnote{Recall that here blocks of a permutation are sets of \emph{indices}.}
of $\gamma$ such that:
\vspace{-6.5pt}
\begin{itemize}\setlength{\itemsep}{-3pt}
	\item  if blocks $\alpha(i)$ and $\alpha(j)$ are not empty, and $i<j$,
then $\alpha(i)$ consists of smaller indices than $\alpha(j)$;

\item as a word, $\alpha(1) \ldots \alpha(n)$ is a factorization of the word
$1\ldots |\gamma|$ (which may include empty factors). 

\item  denoting $\gamma_I$ the pattern corresponding to 
$\gamma_{i_1} \ldots \gamma_{i_{\ell}}$ for any
block $I$ of indices from $i_1$ to $i_{\ell}$ in increasing order, we
have $\pi\{\gamma_{\alpha(1)},\dots,\gamma_{\alpha(n)}\}=\gamma$. 
\end{itemize}
\vspace{-5.5pt}
There are $11$ embeddings of $\gamma = 5\,4\,6\,3\,1\,2$ into
$\pi = 3\,1\,4\,2$, which correspond for instance to the generalized substitutions $\pi\{3241,12,0,0\}$,
$\pi\{3241,0,0,12\}$ and $\pi\{0,0,3241,12\}$ for the
same expression of $\gamma$ as the substitution ${21[3241,12]}$, or $
\pi\{3241,1,0,1\}$ which is the only one corresponding to
$312[3241,1,1]$. Notice that this definition of embeddings conveys the
same notion than in \cite{AA05}, but it is formally different and it will turn to be more
adapted to the definition of the sets $E_{i,j}$.

Equations~\eqref{eqn:const1} to~\eqref{eqn:const3} can be viewed as
Equations~\eqref{eqn:Wc1} to~\eqref{eqn:Wc3} ``decorated'' with
pattern avoidance constraints. These constraints apply to every set
$\pi[\wc_1, \dots, \wc_n]$ that appears in a disjoint union on the
right-hand side of an equation. For each such set, the pattern
avoidance constraints can be expressed by pushing constraints into the
subtrees, using embeddings of excluded patterns in the root $\pi$.  For
instance, assume that $\gamma= 5\,4\,6\,3\,1\,2 \in B^\star$ and $\cal
S_{\cal C}=\{3142\}$, and consider $3142[\wc, \wc,\wc,
\wc]\langle \gamma \rangle$. 
The embeddings of~$\gamma$ in~$3142$ indicates how pattern $\gamma$ can be
found in the subtrees in $3142[\wc, \wc,\wc,
\wc]$. As example the last embedding of the previous example
tells that $\gamma$ can spread over all the subtrees of~$3142$ except
the third. In order to avoid this particular embedding of~$\gamma$, it
is enough to avoid one of the induced pattern $\gamma_I$ on one of the
subtrees. However, in order to ensure that~$\gamma$ is avoided, the
constraints resulting from all the embeddings must be considered and
merged. More precisely, consider a set~$\pi[\C_1, \dots ,\C_n]
\langle \gamma \rangle$, $\pi$ being a simple
permutation. Let~$\{\alpha_1, \dots, \alpha_\ell\}$ be the set of
embeddings of $\gamma$ in $\pi$, each $\alpha_i$ being associated to a
generalized substitution $\gamma = \pi\{\gamma_{\alpha_i(1)},\dots,
\gamma_{\alpha_i(n)}\}$ where $\gamma_{\alpha_i(k)}$ is embedded in
$\pi_k$. Then the constraints are propagated according to the
following equation:
\begin{equation} \label{eq:propagate}
\pi[\C_1, \dots ,\C_n] \langle \gamma \rangle = 
\textstyle\bigcup_{(k_1, \dots, k_\ell) \in K^\pi_\gamma} \pi[\C_1 \langle E_{1,k_1 \dots k_\ell} \rangle, \dots ,\C_n \langle E_{n,k_1 \dots k_\ell} \rangle]
\end{equation}
where $K^\pi_\gamma =\{(k_1, \dots, k_\ell) \in [1..n]^\ell\ |\ \forall i,\ \gamma_{\alpha_i(k_i)} \neq 0\}$ and
$E_{m,k_1 \dots k_\ell }= \{ \gamma_{\alpha_i(k_i)}\ |\ i \in [1..\ell] \text{ and } k_i=m \}$ is a set containing at least $\gamma$ for $(k_1, \dots, k_\ell) \in K^\pi_\gamma$.
In a tuple $(k_1, \ldots, k_{\ell})$ of $K^\pi_\gamma$, $k_i$
indicates a subtree of $\pi$ where the pattern avoidance constraint
($\gamma_{\alpha_i(k_i)}$ excluded) forbids any occurrence of
$\gamma$ that could result from the embedding $\alpha_i$. The set
$E_{m,k_1 \dots k_\ell }$ represents the pattern avoidance constraints
that have been pushed into the $m$-th subtree of $\pi$ by embeddings
$\alpha_i$ of $\gamma$ in $\pi$ where the block $\alpha_i(k_i)$ of
$\gamma$ is embedded into $\pi_m$.

Starting from a finite basis of patterns $B$, Algorithm~\ref{alg:sys-ambigu} describes the whole process to compute an ambiguous system defining the class~$\C = Av(B)$ knowing its set of simple permutations $\s_\C$. 
The propagation of the constraints expressed by Equation~\eqref{eq:propagate} is performed by the
procedure~\textsc{AddConstraints}. It is applied to every set of
the form $\pi[\C_1, \dots , \C_n] \langle B' \rangle$ that appears in the equation defining some
$\wc^{\varepsilon}\langle B' \rangle$ by the procedure~\textsc{ComputeEqn}. Finally,
Algorithm~\ref{alg:sys-ambigu} computes an ambiguous system for a
permutation class $Av(B)$ containing a finite number of simple
permutations: it starts from Equations~\eqref{eqn:const1} to~\eqref{eqn:const3},
and adds new equations to this system calling
procedure~\textsc{ComputeEqn}, until every $\pi[\C_1, \dots ,
\C_n] \langle B' \rangle$ is replaced by some $\pi[\C'_1, \dots ,
\C'_n]$ and until every $\C'_i = \wc^{\varepsilon}\langle B'_i \rangle$ is defined by an equation of the system. All the
sets $B'$ are sets of 
patterns of some permutations in $B$. Since there is only a
finite number of patterns of elements of $B$, there is a
finite number of possible $B'$, and Algorithm~\ref{alg:sys-ambigu} terminates.

\SetKwBlock{PGfunc}{\textsc{AddConstraints}}{end} 
\SetKwBlock{Efunc}{\textsc{ComputeEqn}}{end} 
\begin{algorithm}[h!] 

\KwData{ $B$ is a finite basis of patterns defining
  ${\mathcal C}=Av(B)$ such that 
$\mathcal{S}_\C$ is known and finite.}  \KwResult{A
    system of equations of the form $\D = \bigcup \pi[\D_1, \dots,
    \D_n]$ defining $\C$.}
\Begin{ 
	$\mathcal{E}\leftarrow$ \textsc{ComputeEqn}($(\wc,\Bstar)$) $\cup$ \textsc{ComputeEqn}($(\wc^+,\Bstar)$) $\cup$ \textsc{ComputeEqn}($(\wc^-,\Bstar)$)\\
    \While{ there is a right-only $\wc^{\varepsilon}\langle B'\rangle$ in some equation of $\Sys$} { $\mathcal{E}
      \leftarrow \mathcal{E}~\cup$ \textsc{ComputeEqn}($\wc^{\varepsilon}$, $B'$)
    } }
\bigskip

/* Returns an equation defining $\wc^{\varepsilon}\langle B'\rangle$ as a union of $\pi[\C_1, \dots , \C_n]$ */\\
/* $B'$ is a set of permutations, $\wc^{\varepsilon}$ is given by $\s_\wc$ and $\varepsilon \in \{~~ , +, -\}$ */\\
\Efunc(\params{$\wc^\varepsilon,B'$}){
$\mathcal{E} \leftarrow$ Equation \eqref{eqn:const1} or \eqref{eqn:const2} or \eqref{eqn:const3} (depending on $\varepsilon$) written with $B'$ instead of $\Bstar$\\
\ForEach{ $t=\pi[\C_1, \dots , \C_n] \langle B' \rangle$ that appears in $\mathcal{E}$}
{$t\leftarrow$ \textsc{AddConstraints}$(\pi[\C_1, \dots , \C_n], B' )$} 
\Return $\mathcal{E}$\\
}

\bigskip
/* Returns a rewriting of $\pi[\C_1 \dots \C_n] \langle E \rangle$ as a union $\bigcup \pi[\D_1, \dots \D_n]$ */\\
\PGfunc(\params{$(\pi[\C_1 \ldots \C_n] , E)$}){ 
\lIf{$E = \emptyset$}{return $\pi[\C_1 \dots \C_n]$}\;
\Else{ choose $\gamma \in E$ and compute all the embeddings of $\gamma$ in $\pi$\\ compute $K^\pi_\gamma$ and sets $E_{m,k_1 \dots k_\ell }$ defined in Equation~\eqref{eq:propagate} \\ return $\bigcup_{(k_1, \dots, k_\ell) \in K^\pi_\gamma} \textsc{AddConstraints}(\pi[\C_1 \langle E_{1,k_1 \dots k_\ell} \rangle, \dots ,\C_n \langle E_{n,k_1 \dots k_\ell} \rangle],  E \setminus \gamma)$.}
}

\caption{\textsc{AmbiguousSystem}($B$)\label{alg:sys-ambigu}} 
\end{algorithm}

Consider for instance the class $\C =Av(B)$ for $B=\{1243,2413,531642,41352\}$: $\C$ contains only one simple permutation (namely $3142$), and $\Bstar = \{1243\}$. Applying Algorithm~\ref{alg:sys-ambigu} to this class $\C$ gives the following system of equations: 
\begin{small}
\begin{eqnarray}
\wc\langle1 2 4 3 \rangle &=& 1 \ \cup \ 12[\wc^{+}\langle 1 2 \rangle , \wc\langle1 3 2 \rangle] \ \cup \ 12[\wc^{+}\langle 1 2 4 3 \rangle , \wc\langle2 1 \rangle] \ \cup \ 21[\wc^{-}\langle 1 2 4 3 \rangle , \wc\langle1 2 4 3 \rangle] \nonumber\\
&\ \cup \ & 3 1 4 2[\wc\langle1 2 4 3 \rangle , \wc\langle1 2 \rangle , \wc\langle2 1 \rangle , \wc\langle1 3 2 \rangle] \ \cup \ 3 1 4 2[\wc\langle1 2 \rangle , \wc\langle1 2 \rangle , \wc\langle1 3 2 \rangle , \wc\langle1 3 2 \rangle] \label{eqn:ambigu1}\\
\wc\langle1 2 \rangle &=& 1 \ \cup \ 21[\wc^{-}\langle 1 2 \rangle , \wc\langle1 2 \rangle]\label{eqn:ambigu2}\\
\wc\langle1 3 2 \rangle &=& 1 \ \cup \ 12[\wc^{+}\langle 1 3 2 \rangle , \wc\langle2 1 \rangle] \ \cup \ 21[\wc^{-}\langle 1 3 2 \rangle , \wc\langle1 3 2 \rangle]\label{eqn:ambigu3}\\
\wc\langle2 1 \rangle &=& 1 \ \cup \ 1 2 [\wc^{+}\langle 2 1 \rangle , \wc\langle2 1 \rangle].\label{eqn:ambigu4}
\end{eqnarray}
\end{small}

\vspace*{-2em}

\section{Disambiguation of the system}
\label{sec:disambiguation}

In the above, Equation~\eqref{eqn:ambigu1} gives an ambiguous
description of the class $\wc\langle1 2 4 3 \rangle$. As noticed in
\cite{AA05}, we can derive an unambiguous equation using the
inclusion-exclusion principle: 
{\small $
    \wc\langle1 2 4 3 \rangle = 1 \ \cup \ 1 2 [\wc^{+}\langle 1 2 \rangle , \wc\langle1 3 2 \rangle] \ \cup \ 1 2[\wc^{+}\langle 1 2 4 3 \rangle , 
\wc\langle2 1 \rangle] \ \setminus \ 1 2 [\wc^{+}\langle 1 2 \rangle , \wc\langle2 1 \rangle]  \ \cup \ 2 1 [\wc^{-}\langle 1 2 4 3 \rangle , \wc\langle1 2 4 3 \rangle] \ \cup \ $
$
    3 1 4 2 [\wc\langle1 2 \rangle , \wc\langle1 2 \rangle ,
    \wc\langle1 3 2 \rangle , \wc\langle1 3 2 \rangle] \ \cup \ 3 1 4
    2 [\wc\langle1 2 4 3 \rangle , \wc\langle1 2 \rangle , \wc\langle2
    1 \rangle , \wc\langle1 3 2 \rangle] \ \setminus \ 3 1 4
    2[\wc\langle1 2 \rangle , \wc\langle1 2 \rangle , \wc\langle2 1
    \rangle , \wc\langle1 3 2 \rangle]$}.
The system so obtained contains negative terms in
general. This still gives a system of equations allowing to compute
the generating function of the class. However, this cannot be easily
used for random generation, as the subtraction of combinatorial
objects is not handled by random samplers.
In this section we disambiguate this system to obtain a new positive one: 
the key idea is to replace the negative terms by \emph{complement sets}, hereby 
transforming pattern avoidance constraints into pattern \emph{containment} constraints.

\subsection{General framework}
The starting point of the disambiguation is to rewrite ambiguous
terms like $A \cup B \cup C$ as a disjoint union $(A \cap B
\cap C) \uplus (\bar{A} \cap B \cap C) \uplus (\bar{A} \cap \bar{B}
\cap C) \uplus (\bar{A} \cap B \cap \bar{C}) \uplus (A \cap \bar{B}
\cap C) \uplus (A \cap \bar{B} \cap \bar{C}) \uplus (A \cap B \cap
\bar{C})\textrm{.}$
By disambiguating the union $A \cup B \cup C$ using complement sets
instead of negative terms, we obtain an unambiguous description of 
the union with only positive terms.
But when taking the complement of a set defined by pattern avoidance 
constraints, these are transformed into pattern \emph{containment} constraints. 

Therefore,
for any set $\mathcal{P}$ of permutations, we define the {\em \rs}
$\mathcal{P}\langle E \rangle(A)$ of $\mathcal{P}$ as the set of
permutations that belong to $\mathcal{P}$ and that avoid every pattern
of $E$ and contain every pattern of $A$. This notation will be used
when $\mathcal{P} = \wc^{\varepsilon}$, for $\varepsilon \in \{~~ , +,
-\}$ and $\C$ a permutation class. With this notation, notice also
that for $A=\emptyset$, $\C\langle E \rangle = \C\langle E
\rangle(\emptyset)$ is a standard permutation class.  Restrictions
have the nice feature of being stable by intersection as
$\mathcal{P}\langle E \rangle(A) \cap \mathcal{P}\langle E'
\rangle(A') = \mathcal{P}\langle E \cup E' \rangle(A \cup
A')$. We also define a {\em \rt} to be a set of permutations
described as $\pi[{\mathcal S}_{1},{\mathcal S}_{2},\ldots,{\mathcal S}_{n}]$ where $\pi$ is a simple
permutation or $12$ or $21$ and the ${\mathcal S}_{i}$ are \rss. By uniqueness of the substitution decomposition of a permutation, restriction terms are stable by intersection as well and the intersection is performed componentwise for terms sharing the same root: $\pi[{\mathcal S}_{1},{\mathcal S}_{2},\ldots,{\mathcal S}_{n}] \cap \pi[{\mathcal T}_{1},{\mathcal T}_{2},\ldots,{\mathcal T}_{n}] = \pi[{\mathcal S}_{1}\cap {\mathcal T}_{1},{\mathcal S}_{2}\cap {\mathcal T}_{2},\ldots,{\mathcal S}_{n}\cap {\mathcal T}_{n}]$.

\subsection{Disambiguate\label{sec:disambiguate}}

The disambiguation of the system obtained by
Algorithm~\ref{alg:sys-ambigu} is performed by
Algorithm~\ref{alg:disambiguise}.  It consists in two main operations. One is
the disambiguation of an equation according to the root of the terms that induce ambiguity, which may introduce right-only restrictions. This leads to the second
procedure which computes new equations (that are added to the system) to describe these new restrictions (Algorithm~\ref{alg:ComputeEquationForRestriction}).

As stated in Section~\ref{sec:ambiguous}, every equation $F$ of our
system can be written as $t =1 \cup t_{1} \cup t_{2} \cup t_{3} \ldots
\cup t_{k}$ where the $t_{i}$ are \rts and $t$ is a
\rs.
By uniqueness of the substitution decomposition of a permutation, terms
of this union which have different roots $\pi$ are disjoint. Thus for
an equation we only need to disambiguate unions of terms with same
root.

\begin{algorithm}[t]
\KwData{A ambiguous system $\Sys$ of combinatorial equations \hfill /* obtained by Algo.~\ref{alg:sys-ambigu} */}
\KwResult{An unambiguous system of combinatorial equations equivalent to $\Sys$ }
\Begin{
\While{ there is an ambiguous equation $F$ in $\Sys$ }{
	Take $\pi$ a root that appears several times in $F$ in an ambiguous way\\
	Replace the restriction terms of F whose root is $\pi$ by a disjoint union using Eq.~\eqref{eq:DisambiguateRoot} -- \eqref{eq:14}\\
	\While{ there exists a right-only \rs $\wc^{\varepsilon}\langle E\rangle(A)$ in some equation of $\Sys$}{
		$\Sys \longleftarrow \Sys \bigcup$ \textsc{ComputeEqnForRestriction}($\wc^{\varepsilon}$,$E$,$A$). \hfill /* See Algo.~\ref{alg:ComputeEquationForRestriction} */
}
}
\Return $\Sys$
}
\caption{\textsc{DisambiguateSystem}($\Sys$) \label{alg:disambiguise}}
\end{algorithm}
 
For example in Equation~\eqref{eqn:ambigu1}, there are two pairs
of ambiguous terms which are terms with root $3 1 4 2$ and terms with
root $12$.  Every ambiguous union can be
written in the following unambiguous way: 
\begin{equation}\label{eq:DisambiguateRoot}
	\textstyle\bigcup_{i=1}^{k} t_{i}=\textstyle\biguplus_{X
	  \subseteq [1\ldots k], X \not= \emptyset} \bigcap_{i \in X} t_{i}
	\cap \bigcap_{i \in \overline{X}} \overline{t_{i}},	
\end{equation}
where the \emph{complement} $\overline{t_{i}}$ of a \rt $t_{i}$ is
defined as the set of permutations of $\wc$ whose decomposition tree has the
same root than $t_{i}$ but that do not belong to $t_{i}$. 
Equation~\ref{eq:ComplementTerm} below shows that $\overline{t_{i}}$ is not a term in general but can be expressed as a disjoint
union of terms. By distributivity of $\cap$ over $\uplus$, the above expression can therefore be
rewritten as a disjoint union of intersection of terms. Because terms
are stable by intersection, the right-hand side of Equation~\ref{eq:DisambiguateRoot}
is hereby written as a disjoint union of terms.

For instance, consider terms with root $3142$ in Equation~\eqref{eqn:ambigu1}:
$t_{1} = 3 1 4 2 [\wc\langle1 2 \rangle , \wc\langle1 2 \rangle ,
\wc\langle1 3 2 \rangle , \wc\langle1 3 2 \rangle]$ and $t_{2} = 3 1
4 2 [\wc\langle1 2 4 3 \rangle , \wc\langle1 2 \rangle , \wc\langle2
1 \rangle , \wc\langle1 3 2
\rangle]$.  Equation~\eqref{eq:DisambiguateRoot} applied to
$t_1$ and $t_2$ gives an expression of the form 
$\wc\langle1243\rangle = 1 \cup 12[\ldots] \cup 12[\ldots] \cup 21[\ldots] \cup (t_{1} \cap t_{2}) \uplus (t_{1} \cap \overline{t_{2}}) \uplus (\overline{t_{1}} \cap t_{2})\textrm{.}$

To compute the complement of a term $t$, it is enough to write that 
\begin{equation}\label{eq:ComplementTerm}
\overline{t}=\biguplus_{X \subseteq \{1,\ldots,n\}, X \not= \emptyset} \pi[{\mathcal S}'_{1},\ldots,{\mathcal S}'_{n}] \mbox{ where }{\mathcal S}'_{i} = \overline{{\mathcal S}_{i}}\mbox{  if }i \in X\mbox{ and }{\mathcal S}'_{i} = {\mathcal S}_{i}\mbox{  otherwise},
\end{equation}
with the convention that $\overline{{\mathcal S}_{i}} = \wc^{\varepsilon} \setminus {\mathcal S}_{i}$ for ${\mathcal S}_{i} = \wc^{\varepsilon}\langle E\rangle(A)$.
Indeed, by uniqueness of substitution
decomposition, the set of permutations of $\wc$ that do not belong to $t$ but
whose decomposition tree has root~$\pi$ can be written as the union of
terms $u = \pi[{\mathcal S}'_{1},{\mathcal S}'_{2},\ldots,{\mathcal S}'_{n}]$ where ${\mathcal S}'_{i} = {\mathcal S}_{i}$ or
${\mathcal S}'_{i}=\overline{{\mathcal S}_{i}}$ and at least one \rs ${\mathcal S}_{i}$ must be
complemented. For example $\overline{21[{\mathcal S}_{1},{\mathcal S}_{2}]} =
21[{\mathcal S}_{1},\overline{{\mathcal S}_{2}}] \uplus 21[\overline{{\mathcal S}_{1}},{\mathcal S}_{2}]
\uplus 21[\overline{{\mathcal S}_{1}},\overline{{\mathcal S}_{2}}]$.

The complement operation being pushed from \rts down to
\rss, we now compute~$\overline{{\mathcal S}}$, for a given \rs ${\mathcal S} = \wc^{\varepsilon}\langle
E\rangle(A)$, $\overline{{\mathcal S}}$ denoting the
set of permutations of $\wc^{\varepsilon}$ that are not in~${\mathcal S}$. Notice
that 
given a permutation $\sigma$ of $A$, then any permutation $\tau$ of
$\wc^{\varepsilon}\langle \sigma \rangle$ is in $\overline{{\mathcal
    S}}$ because $\tau$ avoids $\sigma$ whereas permutations of
${\mathcal S}$ must contain $\sigma$. Symmetrically, if a
permutation $\sigma$ is in  $E$ then permutations of
$\wc^{\varepsilon}\langle\rangle(\sigma)$ are in $\overline{{\mathcal
    S}}$. It is straightforward to check that $\textstyle
\overline{\wc^{\varepsilon}\langle E \rangle (A)} = \big[
\bigcup_{\sigma \in E} \wc^{\varepsilon}\langle\rangle(\sigma)\big]
\bigcup \big[ \bigcup_{\sigma \in A}
\wc^{\varepsilon}\langle\sigma\rangle()\big]$.  Unfortunately this
expression is ambiguous. Like before we can rewrite it as an
unambiguous union
\begin{equation} \label{eq:14}
\overline{\wc^{\varepsilon}\langle E \rangle (A)}
= \biguplus_{\underset{X \times Y \not=
    \emptyset\times\emptyset}{{X\subseteq A, Y \subseteq E}}}
\wc^{\varepsilon}\langle X\cup\overline{Y} \rangle (Y \cup \overline{X}) \textrm{, where } \overline{X} = A \setminus X \textrm{ and }\overline{Y} = E \setminus Y \textrm{.}
\end{equation}

In our example (Equations~\eqref{eqn:ambigu1} to~\eqref{eqn:ambigu4}),
only trivial complements appear as every \rs is of the form
$\wc\langle \sigma \rangle()$ or $\wc\langle \rangle(\sigma )$ for
which complements are respectively $\wc\langle \rangle(\sigma )$ and
$\wc\langle \sigma \rangle()$.

All together, for any equation of our system, we are able to rewrite it
unambiguously as a disjoint union of \rts. As noticed before,
some new right-only \rss may appear during this process, for example as
the result of the intersection of several \rss or when complementing
\rss. To obtain a complete system we must compute iteratively
equations defining these new \rss using Algorithm~\ref{alg:ComputeEquationForRestriction} described below.

Finally, the terminaison of Algorithm~\ref{alg:disambiguise} is easily proved. 
Indeed, for all the \rss $\wc^{\varepsilon}\langle E\rangle(A)$ that are
considered in the inner loop of Algorithm~\ref{alg:disambiguise}, every
permutation in the sets $E$ and $A$ is a pattern of some element of
the basis $B$ of $\mathcal{C}$. And since $B$ is finite, there is a finite
number of such restrictions.

\subsection{Compute an equation for a \rs}

Let $\wc^{\varepsilon}\langle E \rangle(A)$ be a \rs. Our goal here is to find
a combinatorial specification of this \rs in terms of smaller \rts
(smaller w.r.t. inclusion). 

If $A = \emptyset$, this is exactly the problem addressed in Section~\ref{sec:addConstraint} and solved by pushing down the pattern avoidance constraints in the procedure \textsc{AddConstraints} of Algorithm~\ref{alg:sys-ambigu}.
Algorithm~\ref{alg:ComputeEquationForRestriction} below shows how to propagate also the pattern \emph{containment} constraints induced by $A \neq \emptyset$.

\SetKwBlock{AMfunc}{\textsc{AddMandatory}}{end} 
\begin{algorithm}[H]
\KwData{$\wc^{\varepsilon}, E,A$ with $E,A$ sets of permutations, $\wc^{\varepsilon}$  given by $\s_\wc$ and $\varepsilon \in \{~~ , +, -\}$.}
\KwResult{An equation defining $\wc^{\varepsilon}\langle E \rangle(A)$ as a union of \rts.}
\Begin{
$F \leftarrow$ Equation \eqref{eqn:Wc1} or \eqref{eqn:Wc2} or \eqref{eqn:Wc3} (depending on $\varepsilon$)\\
\ForEach{$\sigma \in E$}{
	/* This step modifies $F$! */\\
	Replace any restriction term $t$ in $F$ by \textsc{AddConstraints}$(t, \{\sigma\})$\hfill /* See Algo.~\ref{alg:sys-ambigu} */\\
}
\ForEach{$\sigma \in A$}{
	/* This step modifies $F$! */\\
	Replace any restriction term $t$ in $F$ by \textsc{AddMandatory}$(t, \sigma)$ \\
}

\Return $F$ \\
}

\bigskip
\AMfunc(\params{$\pi[{\mathcal S}_1, \dots, {\mathcal S}_n],\gamma$}){
\Return a rewriting of $\pi[{\mathcal S}_1, \dots, {\mathcal S}_n] (\gamma)$ as a union of \rts using Equation~\eqref{eq:AddMandatory}.
}
\caption{\textsc{ComputeEqnForRestriction}$(\wc^{\varepsilon},E,A)$ \label{alg:ComputeEquationForRestriction}} 
\end{algorithm}

The pattern \emph{containment} constraints are propagated by \textsc{AddMandatory}, in a very similar fashion to the pattern \emph{avoidance} constraints propagated by \textsc{AddConstraints}. 
To compute $t(\gamma)$ for $\gamma$ a permutation and $t = \pi[{\mathcal S}_1, \dots, {\mathcal S}_n]$ a restriction term, we first  compute all embeddings of $\gamma$ into $\pi$.
In this case, a permutation belongs to $t(\gamma)$ if and only if at
least one embedding is satisfied. 
Hence, any restriction term $t = \pi[{\mathcal S}_1, \dots, {\mathcal S}_n](\gamma)$ rewrites as a (possibly ambiguous) union as follows:
\begin{equation}\label{eq:AddMandatory}
	\textstyle\bigcup_{i=1}^{\ell} \pi[{\mathcal S}_{1}(\gamma_{\alpha_{i}(1)}),{\mathcal S}_{2}(\gamma_{\alpha_{i}(2)}),\ldots,{\mathcal S}_{n}(\gamma_{\alpha_{i}(n)})],	
\end{equation}
where the  $(\alpha_{i})_{i \in \{1, \ldots, \ell \}}$ are all the embeddings of $\gamma$ in $\pi$ and if $\gamma_{\alpha_{i}(j)}=0$, then  ${\mathcal S}_{j}(\gamma_{\alpha_{i}(j)}) = {\mathcal S}_j$.
For instance, for $t = 2413[{\mathcal S}_{1},{\mathcal S}_{2},{\mathcal S}_{3},{\mathcal S}_{4}]$ and $\gamma = 3214$, there are $9$ embeddings of $\gamma$ into $2413$, and the embedding $2413\{321,1,0,0\}$ contributes to the above union with the term $2413[{\mathcal S}_{1}(321),{\mathcal S}_{2}(1),{\mathcal S}_{3},{\mathcal S}_{4}]$. 

Notice that although the unions of Equation~\ref{eq:AddMandatory} may be ambiguous, they will be transformed into disjoint unions by the outer loop of Algorithm~\ref{alg:disambiguise}. 
Finally, the algorithm produces an unambiguous system which is the result of a finite number of iterations of computing equations followed by their disambiguation.

\section{Conclusion}
\label{sec:ccl}

We provide an algorithm to compute a combinatorial specification 
for a permutation class $\C=Av(B)$, when its basis $B$ and the set of its simple permutations are finite and given as input. The complexity of this algorithm is however still to analyse. In particular, we observe a combinatorial explosion of the number of equations in the system obtained, that needs to be quantified.

Combined with existing algorithms, our procedure provides a full algorithmic chain from the basis (when finite) of a permutation class \C to a specification for \C. This procedure may fail to compute its result, when \C contains an infinite number of simple permutations, this condition being tested algorithmically.

This procedure has two natural algorithmic continuations.
First, with the \emph{dictionnary} of \cite{FlSe09}, the
constructors in the specification of $\mathcal{C}$ can be
directly translated into operators on the generating function $C(z)$
of $\mathcal{C}$, turning the specification into a system of (possibly
implicit) equations defining $C(z)$.  
Notice that, using the inclusion-exclusion principle as
in \cite{AA05}, a system defining $C(z)$ could also be obtained
from an \emph{ambiguous} system describing
$\mathcal{C}$.
Second, the specification can be translated directly
into a Boltzmann uniform random sampler of permutations in
$\mathcal{C}$, in the same fashion as the above dictionnary (see
\cite{DuFlLoSc04}). This second translation is possible
only from an  unambiguous system: indeed,
whereas adapted when considering enumeration sequences, the
inclusion-exclusion principle does not apply when working on the
combinatorial objects themselves. 

When generating permutations with a Boltzmann sampler, complexity 
is measured w.r.t. the size
of the permutation produced (and is linear if we allow a small variation on the size of the output
permutation; quadratic otherwise) and not at all w.r.t. the number of equations in
the specification. In our context, this dependency is of course
relevant, and opens a new direction in the study of Boltzmann random
samplers.

With a complete implementation of the algorithmic chain from $B$ to the specification
and the Boltzmann sampler, one should be able to test conjectures on and study
permutation classes. One direction would be to somehow
measure the randomness of permutations in a given class, by comparing
random permutations with random permutations in a class, or
random permutations in two different classes, w.r.t. well-known
statistics on permutations. Another perspective would be to use the
specifications obtained to compute or estimate the growth rates of
permutation classes, to provide improvements on the known bounds on
these growth rates. We could also explore the possible use 
the computed specifications to provide more
efficient algorithms to test membership of a permutation to a
class.

However, a weekness of our procedure that we must acknowledge is that it fails to be completely general. Although the method is generic and algorithmic, the classes that are fully handled by the algorithmic process are those containing a finite number of simple permutations. By \cite{AA05}, such classes have finite basis (which is a restriction we imposed already), but they also have an  \emph{algebraic} generating function. Of course, this is not the case for every permutation class.
We may wonder how restrictive this framework is, depending on which problems are studied.
First, does it often happen that a permutation class contains
finitely many simple permutations? To properly express what \emph{often} means, a probability distribution on permutation classes should be defined, which is a
direction of research yet to be explored.
Second, we may want to describe some problems (maybe like the distribution of some statistics) for which algebraic permutation classes are representative of all permutation classes.

To enlarge the framework of application of our algorithm, we could explore
the possibility of extending it to permutation
classes that contain an infinite number of simple permutations, but
that are finitely described (like the family of oscillations of
\cite{BRV06} for instance). With such an improvement, more classes would enter our framework, but it would be hard to leave the algebraic case.
This is however a promising
direction for the construction of Boltzmann random samplers for such
permutation classes.

\bibliographystyle{abbrvnat}
\begin{small}
\bibliography{biblio}
\end{small}

\end{document}